\documentclass[reqno]{amsart}
\usepackage{amsmath,amsthm,amssymb}
\usepackage{latexsym}
\usepackage{eucal}

\newtheorem{theorem}{Theorem}[section]

\theoremstyle{definition}

\def\cal{\mathcal}
\let\Re=\undefined
\DeclareMathOperator{\Re}{Re}
\let\Im=\undefined
\DeclareMathOperator{\Im}{Im}

\begin{document}
\title[Multidimensional $L^2$ conjecture: a survey ]{Multidimensional $L^2$ conjecture: a survey}

\author{Sergey A. Denisov}

\email{denissov@math.wisc.edu}

\thanks{{\it Schr\"odinger operator, slowly decaying potential, scattering, absolutely continuous spectrum } }

\address{University of Wisconsin-Madison,
Mathematics Department, 480 Lincoln Dr. Madison, WI 53706-1388, USA}

\maketitle \centerline{\it{To the 70-th birthday of N.K. Nikolskii}}
\begin{abstract}
In this survey, we will give a short overview of the recent progress
on the multidimensional $L^2$ conjecture. It can also serve  as a
quick introduction to the subject. Another survey  was recently
written by Oleg Safronov \cite{saf1} and we highly recommend it.
\end{abstract} \vspace{1cm}

\section{Introduction}

The one-dimensional scattering theory for Schr\"odinger and Dirac
operators is fairly well-understood by now with only a few very
difficult problems left. This progress is mostly due to applying the
tools of complex function theory and harmonic analysis.  In
multidimensional situation, very little is known.

The multidimensional $L^2$ conjecture was suggested by Barry Simon
\cite{simon} and it asks the following.\smallskip

{\bf Conjecture.}  Consider
\[
H=-\Delta+V, \quad x\in \mathbb{R}^d
\]
where $V$ is real valued potential which satisfies
\begin{equation}
\int_{\mathbb{R}^d}\frac{V^2(x)}{1+|x|^{d-1}}dx<\infty \label{l2}
\end{equation}
Is it true that $\sigma_{ac}(H)$ contains the positive
half-line and it is of infinite multiplicity there?\smallskip

One might have to require more local regularity for the potential
just to define $H$ correctly \cite{cycon}, e.g., assuming $V\in L^\infty(\mathbb{R}^d)$ is already good enough.

This conjecture was completely solved only for $d=1$ \cite{deift}.
Even for the case $|V(x)|\lesssim (1+|x|)^{-\gamma}, \gamma=(1-)$
nothing is known. Below, we will discuss some cases in which the
progress was made. We will also briefly explain the methods and
suggest some open problems.

\section{Cayley tree}

The material in this section is taken from \cite{d1,d2, simon2}.
Assume that the Cayley tree $\mathbb{B}$ is rooted with the root
(the origin) denoted by $O$, $O$ has two neighbors and other
vertices have one ascendant and two descendants (the actual number
of descendants is not important but it should be the same for all
points $X\neq O$). The set of vertices of the tree is denoted by
$\mathbb{V}(\mathbb{B})$. For an $f\in
\ell^2(\mathbb{V}(\mathbb{B}))$, define the free Laplacian by
$$
({H}_0f)_n=\sum_{{\rm dist}(i,n)=1}f_i, \quad n\in \mathbb{V}(\mathbb{B})
$$
One can show rather easily  that the spectrum of ${H}_0$ is purely
a.c. on $[-2\sqrt 2, 2\sqrt 2]$. Assume now that $V$ is a bounded
potential on $\mathbb{V}(\mathbb{B})$ so that
$$
{H}={H}_0+V
$$
is well-defined. Denote the spectral measure related to delta
function at $O$ by $\sigma_O$; the density of its absolutely continuous part  is
 $\sigma'_O$. Take $w(\lambda)=(4\pi)^{-1}(8-\lambda^2)^{1/2}$ and let
$\rho_O(\lambda)=\sigma'_O(\lambda)w^{-1}(\lambda)$.

Consider also the probability space on the set of nonintersecting paths in
$\mathbb{B}$ that go from the origin to infinity. This space is
constructed by assigning the Bernoulli random variable to each
vertex and the outcome of Bernoulli trial ($0$ or $1$) then
corresponds to whether the path (stemming from the origin) goes
to the ``left" or to the ``right" descendant at the next step. Notice
also that (discarding a set of Lebesgue measure zero) each path is
in one-to-one correspondence with a point on the interval $[0,1]$ by
the binary decomposition of reals. In this way, the ``infinity"
for $\mathbb{B}$ can be identified with $[0,1]$. For any $t\in
[0,1]$, we can then define the function $\phi$ as
$$
\phi(t)=\sum_{n=1}^\infty V^2(x_n)
$$
where the path $\{x_n\}\subset \mathbb{V}(\mathbb{B})$ corresponds
to $t$. This function does not have to be finite at any point $t$
but it is well-defined and is Lebesgue measurable.  See \cite{d1}
for
\begin{theorem}\label{tree}
For any bounded $V$,
\begin{eqnarray*}\label{bs}
\int\limits_{-2\sqrt 2}^{2\sqrt 2} w(\lambda)\log
\rho_O(\lambda)d\lambda &\ge& \log \mathbb{E} \left\{ \exp\left[
-\frac 14 \sum\limits_{n=1}^\infty V^2(x_n)\right]\right\}\\
&=&\log \int\limits_0^1\exp\left(-\frac{\phi(t)}{4}\right)dt
\end{eqnarray*}
where the expectation is taken with respect to all paths
$\left\{x_n\right\}$ and the probability space defined above. In
particular, if the right hand side is finite, then $[-2\sqrt
2,2\sqrt2]\subseteq \sigma_{ac}({H})$.
\end{theorem}
The proof of  the theorem is based on the adjusted form of sum rules
in the spirit of Killip-Simon \cite{killip-simon}. Higher order sum
rules are applied to different classes of potentials in Kupin
\cite{kupin}.

Notice that $\phi$ is always nonnegative, therefore the right hand
side is bounded away from $-\infty$ iff $V\in \ell^2$ with a
positive probability. This is the true multi-dimensional
$L^2$ condition. The simple application of Jensen's inequality then
immediately implies that the estimate
\[
\int \phi(t)dt=\sum\limits_{n=0}^\infty 2^{-n} \sum\limits_{{\rm
dist}(X,O)=n} V^2(X)<\infty
\]
guarantees  $[-2\sqrt 2,2 \sqrt2]\subseteq \sigma_{ac}({H})$. The
last condition is precisely the analogue of (\ref{l2}) for the
Cayley tree. Indeed, the factor $2^{n}$ is the ``area" of the sphere
of radius $n$ in $\mathbb{B}$ and is exactly the counterpart of
$|x|^{d-1}$ in the same formula.

\section{Slowly decaying oscillating potentials}

There are two different methods to handle this case. \smallskip

{\bf 1.  Asymptotics of Green's function for the complex values of
spectral parameter.}

For simplicity, take $d=3$ and assume that $V$ is supported on the
ball of radius $\rho$ around the origin. Consider the resolvent
$R_z=(H-z)^{-1}, z\in \mathbb{C}^+$ and denote its integral kernel
by $G_z(x,y)$. The approach suggested in \cite{den1} requires the
careful analysis of the asymptotical behavior of $G_z(x,y)$ when
$z\in \mathbb{C}^+$, $y\in \mathbb{R}^3$ are fixed and $x\to\infty$
in some direction. To be more precise, we compare $G_z(x,y)$ to the
unperturbed Green's function
$$G_z^0(x,y)=\frac{\exp (ik|x-y|) }{4\pi |x-y|},\quad z=k^2$$ in the following way.
Take any $f(x)\in L^2(\mathbb{R}^3)$ with a compact support and
define $u=R_zf$. As $V$ is compactly supported, we have
\begin{equation}
\begin{array}{cc}
\displaystyle u(x,k)=&\displaystyle \frac{\exp (ikr)}{r}\left(
A(k,\theta)+\bar{o}(1)\right), \\
\displaystyle \frac{\partial u(x,k)}{\partial r}= &ik
\displaystyle \frac{\exp
(ikr)}{r}\left( A(k,\theta)+\bar{o}(1)\right),\\
\displaystyle & r=|x|,\, \theta=\displaystyle \frac{x}{|x|},
|x|\to\infty
 \end{array}
 \quad {\rm (Sommerfeld's\; radiation\; conditions)}
\end{equation}
Let us call $A$ an amplitude. Clearly, its analysis boils down to
computing the asymptotics for $G_z(x,y)$. The amplitude
$A(k,\theta)$ has the following properties:

\begin{itemize}
\item[1.] { $A(k, \theta)$ is a
vector-function analytic in $k\in \{\Im k>0, \Re k>0\}$.}
\item[2.] {The absorption principle holds, i.e. $A(k,\theta)$ has a
continuous extension to the positive half-line.}

\item[3.] {For the boundary value of the resolvent, we have\\
 \phantom \quad \mbox{$\Im
(R^+_{k^2}f,f)=k\|A(k,\theta)\|_{L^2(\Sigma)}^2, k>0$}. Therefore,
\begin{equation}
\sigma'_f(E)=k\pi^{-1}\|A(k,\theta)\|_{L^2(\Sigma)}^2, E=k^2
\label{factor}
\end{equation}
where $\sigma_f(E)$ is the spectral measure of $f$.}
\end{itemize}
The last formula is the crucial one. The key observation made in
\cite{den1} is that the function $\log
\|A(k,\theta)\|_{L^2(\Sigma)}$ is subharmonic in $k\in \{\Im k\geq
0, \Re k>0\}$. Thus, provided that some rough estimates (uniform in
$\rho$) are available for $\|A(k,\theta)\|_{L^2(\Sigma)}$ away from
real axis, one can use the mean-value formula to get
\begin{equation}
\int_{I} \log \sigma'_f(E)dE>C\label{ent}
\end{equation}
for any interval $I\subset\mathbb{R}^+$. Then, as $C$ is
$\rho$--independent, one  can extend this estimate to the class of
potentials that are not necessarily compactly supported. This
requires using the lower-semicontinuity of the entropy
\cite{killip-simon}. The estimate (\ref{ent}) yields
$\sigma'_f(E)>0$ for a.e. $E>0$ and the statement on the a.c.
spectrum easily follows.

The technical part is to obtain the estimates on the amplitude $A$.
This can be done by the perturbation theory technique. The typical
result \cite{den1} one can obtain this way is
\begin{theorem} \label{imrn1}
Let $Q(x)$ be a $C^1(\mathbb{R}^3)$ vector-field in $\mathbb{R}^3$
and
\[
|Q(x)|<\frac{C}{1+|x|^{0.5+\varepsilon}},\, |div\;
Q(x)|<\frac{C}{1+|x|^{0.5+\varepsilon}}, \varepsilon>0
\]
Then, $H=-\Delta+div\; Q$ has an a.c. spectrum that fills
$\mathbb{R}^+$.\label{ac-spectrum}
\end{theorem}
The decay of potential here is nearly optimal but an additional
oscillatory behavior is also needed. If no oscillation is
assumed, then the whole method breaks down. More about that later.
\bigskip

{\bf 2. Method of Laptev-Naboko-Safronov.}

This very elegant approach was suggested in \cite{saf9} and was
later developed in subsequent publications \cite{saf1, saf2, saf3,
saf4, saf5, saf6, saf7, saf8}. We will again give only a sketch of
the idea. Rewriting the operator in the spherical coordinates we
have
\[
H\sim-\frac{d^2}{dr^2} -\frac{B}{r^2}+V(r,\theta)
\]
where $B$ is Laplace-Beltrami on the unit sphere. Now, let us treat this as a
one-dimensional operator with operator-valued potential
\[
Q= -\frac{B}{r^2}+V(r,\theta)
\]
Denote the projection to the first spherical harmonic by $P_{0}$.
The idea of \cite{saf9} is to write $P_0(H-z)^{-1}P_0$ as

\[
-\frac{d^2}{dr^2} +Q(z)
\]
with nonlocal potential $Q(z)$ and later apply the one-dimensional technique to this operator.
 For example, the following result can be obtained this way

\begin{theorem}
Let $d\geq 3$ and $V(x)$ be such that
\begin{itemize}
\item[1.]{$\lim_{|x|\to \infty} V(x)= 0$}
\item[2.]{$V\in L^{d+1}(\mathbb{R}^{d+1})$}
\item[3.] The Fourier transform of $V$ is well-defined around the origin as $L^2_{\rm loc}$
function, i.e.
\[
\int_{|\xi|<\delta} |\hat{V}(\xi)|^2d\xi<\infty
\]
for some $\delta>0$.
\end{itemize}
Then, $\sigma_{ac}(H)=\mathbb{R}^+$.
\end{theorem}
Notice that the second condition on $V$ is satisfied for, e.g.,
$|V(x)|<C(|x|+1)^{-\gamma}, \gamma>d/(d+1)$. However, the third
condition implies that $V$ either decays fast or oscillates. The
substantial problem with this method is that one needs good bounds
on the discrete negative spectrum, e.g. Lieb-Thirring estimates.
However, the needed estimates can be obtained only under rather
strict assumptions on the decay of potential. In the section
\ref{pencil}, we will explain how this difficulty can be overcame.

 \medskip

\section{Nontrivial WKB correction}

The method of asymptotical analysis of Green's function explained in
the previous section can also be used in two different situations:
when $V_\theta$ is short-range and when $V$ is sparse.

In \cite{gp}, the following result was obtained.
\begin{theorem} Let $d=3$ and
\[
|V(x)|+|x|\cdot |\nabla' V(x)|\lesssim (1+|x|)^{-0.5-}
\]
where $\nabla' V=\nabla V\cdot (x/|x|)^\perp$ is the angular part of
the gradient of $V$. Then, $\sigma_{ac}(H)=\mathbb{R}^+$.
\end{theorem}
The method employed is essentially the same as the one used in
\cite{den1} with one exception: the Green's function asymptotics
contains the well-known WKB-type correction:
\begin{equation}
G_{z}(x,y)\sim \frac{1}{4\pi|x-y|}\exp\left(ik|x-y|+
\frac{1}{2ik}\int_0^{|x|} V(\hat{x}s)ds\right)\label{perelman}
\end{equation}
as $|x|\to\infty$ and $\hat{x}=x/|x|$. This correction becomes a
unimodular factor when $k\in \mathbb{R}$ and so the main arguments
of \cite{den1} go through.\smallskip

In the paper \cite{den2}, quite a different situation was
considered. Take the sequence $R_n$ to be very sparse and consider
the potential $V$ supported on the concentric three-dimensional
shells $\Sigma_n$ with radii $R_n$ and width $\sim 1$.

\begin{theorem}
Assume that $|V(x)|< v_n$ if $x\in \Sigma_n$ and $v_n\in
\ell^2(\mathbb{N})$, then $\sigma_{ac}(H)=\mathbb{R}^+$.
\end{theorem}

The method is again based on the calculation of the asymptotics of
Green's function for the fixed complex $k$. This asymptotics
involves new  and nontrivial WKB factor which is not unimodular for
real $k$ but it is sufficiently regular to apply the same technique.

The WKB correction obtained so far in the literature was always a
quite explicit multiplier. That, however, does not seem to be the
case in general. First of all, in spite of many attempts, no results
on the asymptotics of the Green's function was even obtained for
complex $k$ as long as the only condition assumed of $V$ is the slow
decay: $|V(x)|<C(|x|+1)^{-1+}$. The possible explanation is the
following. Consider $G_{k^2}(x,0)$ for real $k$ as a function of
$\hat{x}=x/|x|$ by going to the spherical coordinates. As
$|x|\to\infty$, the contribution to the $L^2$--norm of this function
coming from the higher angular modes (i.e. the modes of order
$|x|^\alpha, \alpha>0.5$) is likely to be more and more pronounced.
If one makes $k$ complex, this phenomenon can hardly disappear but
the waves corresponding to different angular frequencies have
different Lyapunov exponents even in the free case so their
contributions are all mixed up in the Fourier sum for $G_{k^2}(x,0)$
 making establishing any asymptotics nearly impossible. That questions the
 applicability of the method of \cite{den1} and perhaps
the Green function analysis for real $k$ is needed. The analysis for
real $k$ performed in \cite{christ-kiselev, tao}, however, was never
sufficient to handle the optimal case $V\in L^2(\mathbb{R}^+)$ and
so the prognosis for the resolution of the full $L^2$ conjecture is
rather negative.

What is the WKB correction to Green's function for real $k$ if there
is any asymptotics at all? We do not know yet but there is one
special case when the correction to the asymptotics of the evolution
group $e^{itH}$ was computed. That was done in \cite{den3}.

Assume that $d=3$ and $V$ satisfies the following conditions

{\bf Conditions A:}
\[
 |V|< Cr^{-\gamma},
\left|\frac{\partial V}{\partial r}\right|<Cr^{-1-\gamma},
\left|\frac{\partial^2 V}{\partial
r^2}\right|<Cr^{-1-2\gamma},V(x)\in C^2(\mathbb{R}^3), \quad r=|x|
\]
and
\begin{equation}
1> \gamma>1/2
\end{equation}

The standard by now Mourre estimates immediately show that the
spectrum is purely a.c. on the positive half-line. The question is
what is the long-time asymptotics of $e^{itH}$? Well, the answer  is
not easy as we will see and it requires quite a bit of notation.
Let, again, $B$ be the Laplace-Beltrami operator on the unit sphere
$\Sigma$. Consider the following evolution equation:
\begin{equation}\label{evolution}
iky_\tau(\tau,\theta)=\frac{(By)(\tau,\theta)}{\tau^2}+V(\tau,\theta)y(\tau,\theta),\quad
\tau>0
\end{equation}
where $k\in \mathbb{R}\backslash \{0\}$, $V(\tau,\theta)=V(\tau\cdot \theta)$ is the potential
 written in spherical coordinates, and the function $y(\tau,\theta)\in L^2(\Sigma)$ for
any $\tau> 0$. We introduce $U(k,\tau_0,\tau)f$, the solution of (\ref{evolution})
satisfying an initial condition $U(k,\tau_0, \tau_0)f=f$ where
$\tau, \tau_0> 0$ and $f\in L^2(\Sigma)$. For any $f\in L^2(\Sigma)$, denote
\begin{equation}
W(k,\tau)f=U(k,1,\tau)f \label{doublev}
\end{equation}
and consider the following operator
\[
[\cal{E}(t)f](x)=(2it)^{-3/2}\exp\left[i|x|^2/(4t)\right]\cdot
W(|x|/t, |x|)\left[\hat f(|x|/(2t)\theta)\right] (\hat{x})
\]
 $f\in L^2(\mathbb{R}^3)$.

\begin{theorem}\label{main-theorem}
Assume that $V$ satisfies Conditions A. Then, for any $f\in
L^2(\mathbb{R}^3)$, the following limits (modified wave operators) exist
\[
\cal{W}_{\pm}f=\lim_{t\to \pm \infty} \exp(iHt)\cal{E}(t)f
\]
\end{theorem}
If $V$ is short-range then $W$ can be dropped in the definition of
$\cal{E}$ and the statement of the theorem will still be correct. If
$V$ is long-range but the gradient is short-range, then one can show
that $W$ has the standard multiplicative WKB correction in large
$\tau$ asymptotics similar to the one present in (\ref{perelman}).
In general, the $Wf$ factor can not be simplified much and so the
WKB correction happens to be given by very complicated evolution
equation.

One should expect that in the case when $|V(x)|<C(1+|x|)^{-\gamma},
\gamma\in (0.5,1)$ even more complicated evolution equation appears
both in the spatial  asymptotics of the Green's function and in the
long-time asymptotics of $e^{itH}$. Meanwhile, proving this seems to
be a monumental task as the statement like that even in one-dimensional
case holds not for all $k\neq 0$ but rather for Lebesgue a.e. $k$.

\section{Ito's stochastic equation and modified Harmonic measure}
The resolution of $L^2$ conjecture for the Cayley tree suggests that
may be the condition (\ref{l2}) can be relaxed. Although the
discussion above might have somewhat sobering effect on the reader,
one can hope to at least try to exploit the idea of introducing the
right space of paths. One step in this direction was made in
\cite{saf3} where the Laptev-Naboko-Safronov method was adjusted to
the case when potential is small in the cone. Further progress was
made in the paper \cite{den4}. Recall that for the Cayley tree one
can very naturally introduce the probability space of paths escaping
to infinity. Then, one knows that if the potential is small
($\ell^2$ or is just zero) with positive probability, then the a.c.
spectrum is present. What is an analog of this probability space in
the Euclidean case?

Consider the Lipschitz vector field
\begin{equation}\label{eee01}
 p(x)=\left(\frac{I'_\nu(|x|)}{I_\nu(|x|)}-\nu |x|^{-1}\right)\cdot
\frac{x}{|x|}, \quad \nu=(d-2)/2
\end{equation}
where $I_\nu$ denotes the modified Bessel function. Then, fix any
point $x^0\in\mathbb{R}^d$ and define the following stochastic
process
\begin{equation}
dX_t=p(X_t)dt+dB_t, \quad X_0=x^0 \label{stochastic}
\end{equation}
with the drift given by $p$. The solution to this diffusion process
exists and all trajectories are continuous  and escape
to infinity almost surely.

\begin{theorem}\label{the1}
Assume that $V$ is continuous nonnegative bounded function and
\[
\mathbb{E}_{x^0} \left[ \exp\left(- \int\limits_0^\infty
V(X_\tau)d\tau\right)\right]>0
\]
for some $x_0$. Then, $\mathbb{R}^+\subseteq \sigma_{ac}(H)$.
\end{theorem}

The positivity of the expectation means that with positive probability we have
 $V(X_t)\in L^1(\mathbb{R}^+)$. The application of Jensen's inequality immediately yields that
\[
\int_\mathbb{R} \frac{V(x)}{|x|^{d-1}+1}dx<\infty
\]
implies the preservation of the a.c. spectrum. The method used to
prove theorem~\ref{the1} is again more or less rephrasing of the one
from  \cite{d1} but the language is probabilistic.

The question now is how to compute those probabilities. The
reasonable simplification here is to assume that $V$ is supported on
some complicated set (say, a countable collection of balls) and then
study when is the probability to hit this set smaller than one. This problem
was addressed by introducing the suitable
 potential theory and by proving the estimates on the modified
 Harmonic measure. The interesting aspect of this analysis is in
 relating the geometric properties of support of $V$ to the scattering
 properties of the medium. This is too deep and technical a subject to
 try to state here the relevant results so we refer the reader to the original paper.

\section{Hyperbolic Schr\"{o}dinger pencils \label{pencil}}
Consider the Schr\"odinger operator
\[
H_\lambda=-\Delta+\lambda V
\]
with the coupling constant $\lambda$ and decaying $V$. The study of
its resolvent
\[
R_z=(H_\lambda-k^2)^{-1}, \quad z=k^2
\]
is often complicated by the presence of the negative discrete
spectrum for $H_\lambda$. For example, proving the sum rules requires the
Lieb-Thirring bounds.

In \cite{den5}, the following idea was suggested. Instead of
inverting the operator $H-k^2$, let us try to invert
$P(k)=-\Delta+k\mu V-k^2$, where $\mu$ is a fixed constant. In other
words, we make the coupling constant momentum-dependent: $
\lambda=k\mu $. The operator $P(k)$ is a hyperbolic pencil and
$P(k)$ is invertible for all $k\in \mathbb{C}^+$. The function
$P^{-1}(k)$ is analytic there and so one has no problems with poles.
The study of the Green's function for $P^{-1}(k)$ is more
straightforward and this Green function (call it $M_k(x,y,\mu)$)
agrees with the Schr\"odinger Green's function as long as the
potential $V$ is compactly supported:
\[
G_{k^2}(x,y,k\mu)=M_{k}(x,y,\mu)
\]
The Fubini theorem then allows to translate results obtained for the
Schr\"odinger pencil to the results for the original Schr\"odinger
operator. This idea greatly expands the class of potentials that can
be treated but the results hold only for generic coupling constant.
Below we list three theorems that can be obtained this way. The
first two are taken from \cite{den5} and assume $d=3$; the third one
is from \cite{den6}.

\begin{theorem}
Assume
\[
V(x)=div\, Q(x)
\]
where the smooth vector field $Q(x)$ satisfies
\[
Q(x), |D Q(x)|\in L^\infty (\mathbb{R}^3), \quad \int_{\mathbb{R}^3}
\frac{|Q(x)|^2}{|x|^2+1}dx <\infty
\]
Then for a.e. $\lambda$, $\mathbb{R}^+\subseteq
\sigma_{ac}(H_\lambda)$.\label{theorem2}
\end{theorem}
This theorem, in contrast to theorem \ref{imrn1}, does not assume
the pointwise decay of $V$ and the result obtained is sharp in terms
of the decay of $Q$.

\begin{theorem}
Assume $V(x)$ is bounded and
\[
\int\limits_1^\infty r|v(r)|^2dr<\infty
\]
for $v(r)=\sup_{|x|=r} |V(x)|$. Then for a.e. $\lambda$,
$\sigma_{ac}(H_\lambda)=\mathbb{R}^+$.\label{theorem3}
\end{theorem}
This result is interesting in that it covers the case when the
nontrivial WKB correction to the Green's function asymptotics can be
present. However, there is no any need to establish it over specific
direction as only the estimate on the angular average of Green
function
\[
\int_{|x|=r} |M_k(x,0,\mu)|^2d\sigma_x\sim e^{2\Im k r}
\]
is used. In this case, the oscillation of $M_k(x,0,\mu)$ in the
angular variable is rather weak and this is what makes the analysis
possible.

Yet another result can illustrate the power of this technique
\begin{theorem}
Assume that $V$ is continuous bounded function and
\[
\mathbb{E}_{x^0} \left[ \exp\left(- \int\limits_0^\infty
|V(X_\tau)|d\tau\right)\right]>0
\]
for some $x_0$. Then, $\mathbb{R}^+\subseteq \sigma_{ac}(H_\lambda)$
for a.e. $\lambda$.
\end{theorem}
As one can see a rather unnatural requirement of $V$ to be
nonnegative present in the theorem~\ref{the1} is now removed.

In several recent publication, the idea of making the coupling
constant momentum dependent was applied in combination with other
interesting techniques, see e.g. \cite{saf10,saf11}.

{\bf Remark.} It is well-known that the analysis of the
one-dimensional Schr\"{o}dinger operator is more technically
involved than the analysis of, say, Dirac operator or the Krein
system. In fact, the one-dimensional differential equation for the
Schr\"{o}dinger pencil considered above happens to be identical to
the Dirac operator (and automatically to the Krein system). So, it
is not so unexpected that in multidimensional case this trick makes
analysis simpler.

The following question is quite natural in view of the results
listed above.

{\bf Question.} Assume $V$ decays in some way. How does the a.c.
spectrum of $H_\lambda$ as a set depend on the value of $\lambda\neq
0$. That boils down to studying the dependence of
\[
F_\lambda(z)=\Im \langle (H_{\lambda}-z)^{-1}f,f \rangle
\]
on $\lambda$ around the regular points $z\in \mathbb{R}^+$. For the
one-dimensional case, the analysis in \cite{christ-kiselev} reveals
that for $V\in L^p, p<2$ there is a $\lambda$-independent set of
energies of the full Lebesgue measure which supports the a.c.
spectrum of $H_\lambda$.

\section{Possible directions}

The $L^2$ conjecture the way it is stated does not say much about,
say, Schr\"{o}dinger dynamics so it is possible that there are some
soft analysis arguments that can nail it. The good example of the
 soft analysis is the so-called sum rules for
the Jacobi matrices \cite{killip-simon}. More mature approach,
  in our opinion, is to try to control some quantities relevant for scattering: the Green's
  function, evolution group, etc. That, most likely, will require application
  of hard analysis methods. In this section, we will try to explain what kind
   of technical difficulties one stumbles  upon when trying to address these questions.

Some special evolution equations seem to provide an adequate model
for understanding of what is going on and one very important example
of these equations is (for $d=2$)
\begin{equation}\label{evol1}
iu_t(t,\theta,k)=k\frac{\partial^2_{\theta\theta}u(t,\theta,k)}{(t+1)^2}
+V(t,\theta)u(t,\theta,k), \quad u(0,\theta,k)=u_0(\theta)
\end{equation}
where $|V(t,\theta)|\lesssim (1+t)^{-\gamma}, \gamma\in (0.5,1)$ and
$k$ and $V$ are real. This equation appears (check
(\ref{evolution})) in the WKB correction for the dynamics of
$e^{itH}$ and it is likely to be the right correction for the
spatial asymptotics of the Green's function.

 What can be said
about the solution $u(t,\theta,k)$ as $t\to\infty$? The $L^2$-norm
 $\|u\|_{L^2(\mathbb{T})}$ is preserved in time but how about the growth of Sobolev norms? The
conjecture stated in \cite{den9} is that generically in $k$ we
should have
\[
\|u(t,\theta,k)\|_{H^1(\mathbb{T})}\lesssim t, \quad t\to\infty
\]
thus the transfer of the $L^2$ norm to higher modes happens in a
controlled way and that prevents the resonance formation.

Another important quantity to study is how concentrated the function
$u(t,\theta,k)$ can be in the $\theta$ variable. That can be
controlled by quantities like
\[
I_1=\|u(t,\theta,k)\|_{L^p(\mathbb{T})}, p>2;
\]
or
\[
I_2=\int_{\mathbb{T}} \log |u(t,\theta,k)|d\theta
\]
or
\[
I_3(\delta)=\inf_{|\Omega|>\delta} \int_\Omega
|u(t,\theta,k)|^2d\theta
\]
One might guess that for typical $k$ one has: $I_1$ is bounded in
$t$ and/or $I_2>-C$ and/or $I_3(\delta)>C\delta$ as long as $\delta>0$.

What is the technical difficulty in the analysis of (\ref{evol1})?
Notice first that by going on the Fourier side in $\theta$ one gets
the infinite system of ODE's coupled to each other through
$\hat{V}$. The differential operator will become the diagonal one
and the gaps between the eigenvalues will decrease in $t$. This
deterioration of the gaps is the key signature of the
multidimensional case. It necessitates handling increasing number of
frequencies at once and this is the hardest part of the analysis.
One should notice that the evolution equation with only two
frequencies interacting with each other, e.g.

\[
iX_t= \left[
\begin{array}{cc}
0 & V(t)\\
V(t) & k
\end{array}
\right] X,\quad X(0,k)=I
\]
can be handled by Harmonic analysis methods developed in
\cite{christ-kiselev, tao}.  Anyhow, the equation (\ref{evol1}) is
poorly understood and its analysis is very complicated. Some
progress was made in \cite{den9} but there is clearly a long way to
go.\bigskip

{\bf Conclusion.} We hope that this survey makes a good point that
the wave propagation through the medium with slowly decaying
potential is an interesting physical phenomenon with WKB correction
given by evolution equation. The phenomenon of resonances appearing
for some  energies becomes far more complicated in multidimensional
case and the oscillatory behavior of Green's function is just
another manifestation of that.

Based on the literature published in the last five years, it appears
that the number of mathematicians actively working on the problem
does not exceed number three so hopefully this review will attract
more interest to this beautiful subject.

\section{Acknowledgment}

This research was supported by NSF grants DMS-1067413 and DMS-0635607.

\end{document}